\documentclass[11pt, a4paper]{amsart}
\usepackage[top=35truemm,bottom=35truemm,left=35truemm,right=35truemm]{geometry}
\usepackage{setspace} 
\usepackage{amsmath, amsfonts,amsthm,amssymb,mathrsfs}
\usepackage{enumerate}
\usepackage[dvipdfmx]{graphicx}
\usepackage{ascmac}
\usepackage[arrow,matrix]{xy}
\usepackage{color}
\usepackage{array}
\usepackage{pifont}
\usepackage{longtable}
\usepackage{mathtools}

\newtheorem{theo}{Theorem}

\newtheorem{prop}{Proposition}
\newtheorem{cor}{Corollary}

\def\H{\mathbb{H}}
\def\R{\mathbb{R}}
\def\E{\mathbb{E}}
\def\U{\mathbb{U}}

\begin{document}

\title[ ]{On the growth rate of ideal Coxeter groups in hyperbolic 3-space}
\author{Yohei Komori}
\address{Department of Mathematics, School of Education, Waseda University, Nishi-Waseda 1-6-1, Shinjuku, Tokyo 169-8050, Japan}
\email{ykomori@waseda.jp}
\author{Tomoshige Yukita}
\address{Department of Mathematics, School of Education, Waseda University, Nishi-Waseda 1-6-1, Shinjuku, Tokyo 169-8050, Japan}
\email{yshigetomo@suou.waseda.jp}
\date{}
\thanks{}

\begin{abstract}
We study the set $\mathcal{G}$ of growth rates of of ideal Coxeter groups in hyperbolic 3-space which consists of real algebraic integers greater than 1.
We show that (1) $\mathcal{G}$ is unbounded above while it has the minimum,  
(2) any element of $\mathcal{G}$ is a Perron number, and (3) growth rates of of ideal Coxeter groups with $n$ generators
are located in the closed interval $[n-3, n-1]$.
\end{abstract}

\maketitle


\setstretch{1.1}
\section{Introduction}
The upper half space $\H^3=\{ x=(x_1,x_2,x_3)\in \R^3 \,|\, x_3>0\}$ 
with the metric $|dx|/x_3$ is a model of  hyperbolic $3$--space, so called the {\em upper half space model}.
The Euclidean plane $\E^2=\{ x=(x_1,x_2,x_3)\in \R^3 \,|\, x_3=0\}$ and the point at infinity $\infty$ consist of 
the boundary at infinity $\partial \H^3$ of $\H^3$.
A subset $B \subset \H^3$ is called a {\em hyperplane} of $\H^3$ if it is a Euclidean hemisphere or a  half plane orthogonal to $\E^2$. 
When  we restrict the hyperbolic metric $|dx|/x_3$ of $\H^3$ to $B$, it becomes  a model of  hyperbolic plane.
We define a {\em polytope} as
a closed domain $P$ of $\H^3$ 
which can be written as the intersection of finitely many closed half spaces $H_B$ bounded by hyperplanes $B$, say
$P=\bigcap H_B$.
In this presentation of $P$, we assume that $F_B=P \cap B$ is a hyperbolic polygon of $B$.
$F_B$ is called a {\em facet} of $P$, and $B$ is called the {\em supporting hyperplane} of $F_B$.
If the intersection of  two facets $F_{B_1}$ and $F_{B_2}$ of  $P$ consists of a geodesic segment, it is called an {\em edge} of $P$,
while if the intersection $\bigcap F_B$ of more than two facets is a point,  then it is called a {\em vertex} of $P$.
If $F_{B_1}$ and $F_{B_2}$ intersect only at a point of the boundary $\partial \H^3$ of $\H^3$,
it is called an {\em ideal vertex} of $P$.
A  polytope $P$ is called {\em ideal} if all of its vertices are ideal.
Related to Jakob Steiner's problem on the combinatorial characterization of polytopes inscribed in the two-sphere $S^2$,
ideal  polytopes in hyperbolic 3-space $\mathbb{H}^3$ has been studied extensively  \cite{G, R}.

A {\em horosphere} $\Sigma$ of $\H^3$ based at $v \in \partial \H^3$ is defined by a Euclidean sphere in $\H^3$ tangent to $\E^2$ at $v$
when $v \in \E^2$, or a Euclidean plane in $\H^3$ parallel to $\E^2$ when $v=\infty$.
When we restrict the hyperbolic metric of $\H^3$ to $\Sigma$, it becomes  a model of  Euclidean plane.
Let $v\in \partial \H^3$ be an ideal vertex of a  polytope $P$ in $\H^3$ and $\Sigma$ be a horosphere of $\H^3$ based at $v$ 
such that $\Sigma$ meets just the facets of $P$ incident with $v$. 
Then the {\em vertex  link} $L(v) := P \cap \Sigma$ of $v$ in $P$ is a  Euclidean convex polygon in the horosphere $\Sigma$. 
If $F_{B_1}$ and $F_{B_2}$ are adjacent facets of $P$ incident with $v$, then the Euclidean dihedral angle between $F_{B_1} \cap \Sigma$ and $F_{B_2} \cap \Sigma$
in $\Sigma$  is equal to the hyperbolic dihedral angle between the supporting hyperplanes ${B_1}$ and ${B_2}$ in $\H^3$
 (cf. \cite[Theorem 6.4.5]{R}).

In this paper we consider ideal {\em Coxeter} polytopes $P$,  
which mean that the dihedral angles of edges of $P$ are submultiples of $\pi$.
Since any Euclidean Coxeter polygon is a rectangle or a triangle with dihedral angles $(\pi/2, \pi/3, \pi/6), (\pi/2, \pi/4, \pi/4) $ or $ (\pi/3, \pi/3, \pi/3)$, 
we see that 
the dihedral angles of an ideal Coxeter polytope must be $\pi/2, \pi/3, \pi/4$ or  $\pi/6$.

Any Coxeter polytope $P$  is a fundamental domain of the discrete group $\Gamma$ generated by the set $S$ 
consisting of the reflections with respects to its facets.
We call $(\Gamma, S)$ the {\em Coxeter system} associated to $P$.
In this situation we can define the {\em word length $\ell_S(x)$ of $x \in \Gamma$ with respect to $S$}
by the smallest integer $k \geq 0$ for which there exist $s_1, s_2, \cdots, s_k \in S$ such that $x=s_1s_2 \cdots s_k$.
The {\em growth function} $f_S(t)$ of $(\Gamma, S)$ is the formal power series $\sum_{k=0}^{\infty} a_kt^k$
where $a_k$ is the number of elements $g \in \Gamma$ satisfying $\ell_S(g)=k$.
It is known that 
the {\em growth rate} of $(\Gamma, S)$, $\tau =\limsup_{k \rightarrow \infty} \sqrt[k]{a_k}$ is bigger than $1$  (\cite{dlH1}) and  less than or equal to  the cardinality $|S|$ of $S$
from the definition.  
By means of Cauchy-Hadamard formula, 
the radius of convergence $R$ of  $f_S(t)$ is the reciprocal of $\tau$, i.e. $1/|S| \leq R < 1$.
In practice the growth function $f_S(t)$ which is analytic on $|t|<R$ extends to a rational function $P(t)/Q(t)$  on $\mathbb{C}$
by analytic continuation
where $P(t), Q(t) \in \mathbb{Z}[t]$ are relatively prime.
There are formulas due to Solomon and Steinberg 
to calculate the rational function $P(t)/Q(t)$ from the data of finite Coxeter subgroups of  $(\Gamma, S)$ (\cite{So, St}. See also~\cite{He}).

\begin{theo}{\rm (Solomon's formula)}\\
The growth function $f_S(t)$ of an irreducible finite Coxeter group $(\Gamma, S)$ can be written as 
$f_S(t)=[m_1+1,m_2+1, \cdots, m_k+1]$ where $[n] =1+t+ \cdots +t^{n-1}, [m,n]=[m][n]$,etc. and $\{m_1, m_2, \cdots, m_k \}$
is the set of exponents of $(\Gamma, S)$.
\end{theo}

\begin{theo}{\rm (Steinberg's formula)}\\
Let $(\Gamma, S)$ be a hyperbolic Coxeter group.
Let us denote the Coxeter subgroup of $(\Gamma, S)$ generated by the subset $T\subseteq S$ by $(\Gamma_T,T)$, and denote its growth function by $f_T(t)$.
Set $\mathcal{F}=\{T\subseteq S \;:\; \Gamma_T$ is finite $\}$. Then
$$
\frac{1}{f_S(t^{-1})}=\sum _{T \in \mathcal{F}} \frac{(-1)^{|T|}}{f_T(t)}.
$$
\end{theo}

In this case, $t=R$ is a pole of  $f_S(t)=P(t)/Q(t)$.
Hence $R$ is a real zero of the denominator $Q(t)$ closest to the origin $0 \in \mathbb{C}$
of all zeros of $Q(t)$.
Solomon's formula implies that $P(0)=1$. Hence $a_0=1$ means that $Q(0)=1$.
Therefore $\tau>1$, the reciprocal of $R$, becomes a real algebraic integer
whose conjugates have moduli less than or equal to the modulus of $\tau$. 
If  $t=R$ is the unique zero of $Q(t)$ with the smallest modulus, then $\tau>1$ is a real algebraic integer whose conjugates have moduli less than 
the modulus of $\tau$:
such a real algebraic integer is called a {\em Perron number}.

The following result is a criterion for growth rates to be Perron numbers.

\begin{prop}\rm{(\cite{KU}, Lemma1)}\\
\label{prop:KU}
 Consider the following polynomial of degree $n \geq 2$
 $$
 g(t)=\sum _{k=1}^n a_k t^k -1,
 $$
 where  $a_k$ is a non-negative integer.
 We also assume that the greatest common divisor of $\{k \in \mathbb{N} \; |\; a_k \neq 0 \}$ is $1$.
 Then there is a real number $r_0$,  $0<r_0<1$ which is the unique zero of $g(t)$ having the smallest absolute value of all zeros of $g(t)$.
\end{prop}

In this paper we study the set $\mathcal{G}$ of growth rates of of ideal Coxeter polytopes in hyperbolic 3-space.
We will show that (1) $\mathcal{G}$ is unbounded above while it has the minimum,  
(2) any element of $\mathcal{G}$ is a Perron number, and (3) growth rates of ideal Coxeter polytopes with $n$ facets
are located in the closed interval $[n-3, n-1]$.

\section{Ideal Coxeter polytopes with 4 or 5 facets in $\mathbb{H}^3$ }

Let $p, q, r$ and $s$ be the number of edges with dihedral angles $\pi/2,\; \pi/3,\; \pi/4$, and $\pi/6$ of an ideal Coxeter polytope $P$  in $\mathbb{H}^3$.
By Andreev theorem \cite{A}, we can classify ideal Coxeter polytopes with $4$ or $5$ facets, and calculate the growth functions $f_S(t)$ of $P$
by means of Steinberg's formula and also growth rates,  see Table 1.
Every denominator polynomial has a form $(t-1)H(t)$ and all coefficients of $H(t)$ satisfies the condition of Proposition \ref{prop:KU}, so that the growth rates of 
ideal Coxeter polytopes with $4$ or $5$ facets are Perron numbers.

\begin{table}[htbp]
\begin{flushleft}
\begin{tabular}{|c|c|c|c|}
\hline
$(p, q, r, s)$ & combinatorial type & denominator polynomial & growth rate\\
\hline
$(2, 2, 0, 2)$ & simplex & $(t-1) (3 t^5+t^4+t^3+t^2+t-1)$ &  2.03074 \\
\hline
$(2, 0, 4, 0)$ & simplex & $(t-1) (3 t^3+t^2+t-1)$ & 2.13040 \\
\hline
$(0, 6, 0, 0)$ & simplex & $(t-1) (3 t^2+t-1)$ & 2.30277 \\ 
\hline
$(4, 2, 0, 2)$ & pyramid &$(t-1) (4 t^5+t^4+2 t^3+t^2+2 t-1)$ & 2.74738 \\ 
\hline
$(4, 0, 4, 0)$ & pyramid  & $(t-1) (4 t^3+t^2+2 t-1)$ & 2.84547 \\ 
\hline 
$(2, 5, 0, 2)$ & prism & $(t-1) (5 t^5+2 t^4+t^3+3 t^2+2 t-1)$ & 3.16205 \\ 
\hline
\end{tabular}
\end{flushleft}
\caption{}
\end{table}

\begin{figure}[htbp]
\begin{center}
 \includegraphics [width=170pt, clip]{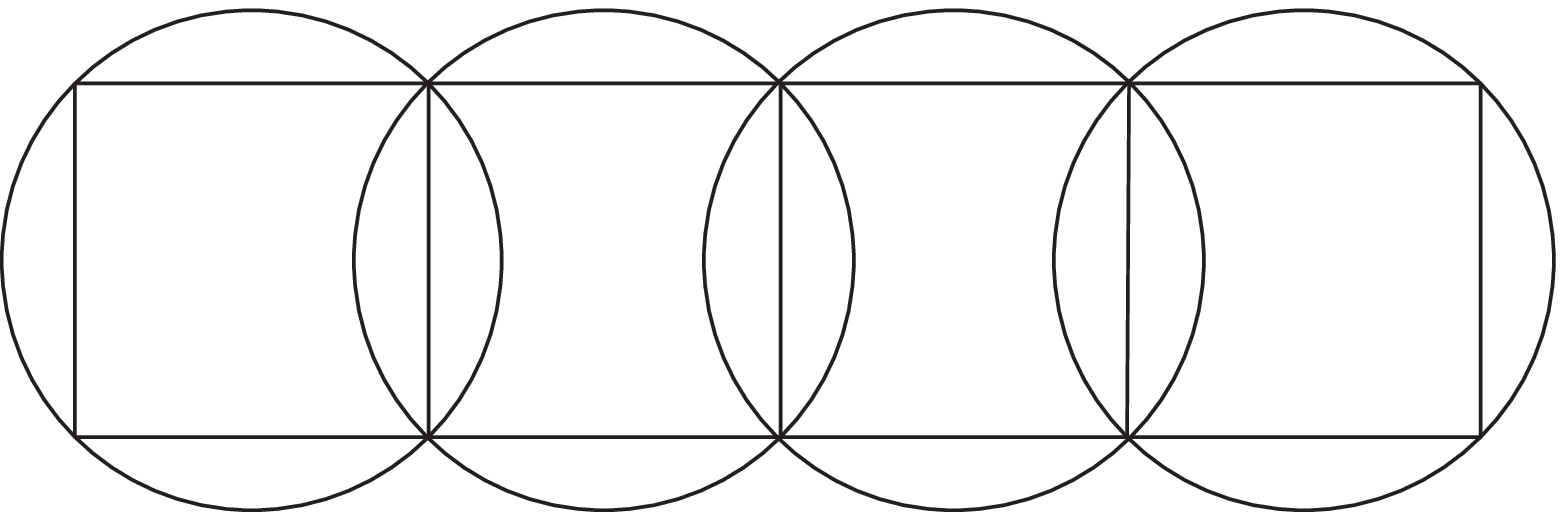}
\end{center}
\caption{}
\label{upperhalf}
\end{figure} 

After glueing $m$ copies of  the  ideal Coxeter pyramid with $p=r=4$ along their sides successively,
we can construct a hyperbolic ideal Coxeter polytope $P$ with $n=m+4$ facets.
In Figure 1 we are looking at   the ideal Coxeter polytope $P$ with $8$ facets
from the point at infinity $\infty$,
which consists of 4 copies of ideal Coxeter pyramid with $p=r=4$ whose apexes are located at  $\infty$;
squares represent bases of pyramids and disks are supporting hyperplanes of these bases.
The growth function of $P$  has the following denominator polynomial
$$
 (t-1)H(t)=(t-1)(2 (n - 3) t^3 + (n - 4) t^2 + (n - 3) t - 1).
$$
Then all coefficients of $H(t)$ except its constant term 
are non-negative. 
Therefore we can apply Proposition \ref{prop:KU} to conclude  that the growth rate of $P$ is a Perron number.
Moreover  $H(t)$ has a unique zero on the unit interval $[0,1]$ and the following inequalities hold:
$$
H(\frac{1}{n-3})=\frac{n-2}{(n-3)^2}>0, \;\; H(\frac{1}{n-1})=\frac{-n^2+n-4}{(n-1)^3}<0.
$$
They imply that the growth rate $\tau$ of $P$  with $n$ facets satisfies
$$
n-3 \leqq \tau \leqq n-1,
$$ 
which will be generalized in the next section. In particular
\begin{prop}
The set $\mathcal{G}$ of growth rates of 3-dimensional hyperbolic ideal Coxeter polytopes is unbounded above.
\end{prop}

\section{The growth rates of ideal Coxeter polytopes in $\mathbb{H}^3$}

Recall that  $p, q, r$ and $s$ be the number of edges with dihedral angles $\pi/2,\; \pi/3,\; \pi/4$, and $\pi/6$ of an ideal Coxeter polytope $P$  in $\mathbb{H}^3$.
By means of Steinberg's formula, we can calculate the growth function $f_S(t)$ of $P$ as
$$
1/f_S(1/t)=1-n/[2]+p/[2,2]+q/[2,3]+r/[2,4]+s/[2,6],
$$
where $[2,3]=[2][3]$, etc.
It can be rewritten as
\begin{eqnarray*}
1/f_S(t) &=& 1-nt/[2]+pt^2/[2,2]+qt^3/[2,3]+rt^4/[2,4]+st^6/[2,6]\\
& = & \frac{1}{[2,2,3,4,6]}G(t),
\end{eqnarray*}
where
$$
G(t)=[2,2,3,4,6]-nt[2,3,4,6]+pt^2[3,4,6]+qt^3[2,4,6]+rt^4[2,3,6]+st^6[2,3,4].
$$

\begin{prop}\label{prop:serre}
Put $a=p/2, \; b=q/3, \; c=r/4, \; d=s/6$. Then
\begin{equation} \label{eq: serre}
a+b+c+d=n-2.
\end{equation}
\end{prop}

\proof
By a result of Serre (\cite{Se}. See also \cite{He})
$$
G(1)=[2,3,4,6](1)(2-n+p/2+q/3+r/4+s/6)=0
$$
\qed

By using this equality (\ref{eq: serre}) we represent $H(t)=G(t)/(t-1)$ as
\begin{eqnarray*}
H(t) &=& -[2,3,4,6]+at[3,4,6]+bt(2t+1)[2,4,6]\\
&&+ct(3t^2+2t+1)[2,3,6]+dt(5t^4+4t^3+3t^2+2t+1)[2,3,4]\\
&=&-1 + (-4 + a + b + c + d) t + (-9 + 3 a + 5 b + 5 c + 5 d) t^2 \\
&&+ (-15 + 6 a + 11 b + 14 c + 14 d) t^3 + (-20 + 9 a + 17 b + 25 c + 29 d) t^4 \\
&&+ (-23 + 11 a + 22 b + 33 c + 49 d) t^5 + (-23 + 12 a + 24 b + 36 c + 66 d) t^6 \\
&&+ (-20 + 11 a + 23 b + 35 c + 71 d) t^7 + (-15 + 9 a + 19 b + 31 c + 61 d) t^8 \\
&&+ (-9 + 6 a + 13 b + 22 c + 40 d) t^9 + (-4 + 3 a + 7 b + 11 c + 19 d) t^{10}\\
&& + (-1 + a + 2 b + 3 c + 5 d) t^{11}.
\end{eqnarray*}

From this formula we have the following result (see also \cite{N}, Theorem 3).

\begin{theo}
The growth rates of ideal Coxeter polytopes in $\mathbb{H}^3$ are Perron numbers.
\end{theo}

\proof
When $n$ the number  of facets satisfies 
$n \geqq 6$, the equality (\ref{eq: serre}) of Proposition
\ref{prop:serre} implies $a+b+c+d=n-2 \geqq 4$. Then all coefficients of $H(t)$ except its constant term 
are non-negative. Hence Proposition \ref{prop:KU} implies the assertion.
For $n=4,5$, this claim was already proved in the previous section. 
\qed

Moreover the equality (\ref{eq: serre}) induces the following two functions
$H_1(t)$ and $H_2(t)$ satisfying
$H_1(t) \leqq H(t) \leqq H_2(t)$ for any $t >0$:

\begin{eqnarray*}
H_1(t) & = &-1 + (-4 +( n - 2) t + (-9 + 3 (n - 2)) t^2 + (-15 + 6 (n - 2)) t^3 \\
&&+ (-20 + 9 (n - 2)) t^4 + (-23 + 11 (n - 2)) t^5 + (-23 + 12 (n - 2)) t^6 \\
&&+ (-20 + 11 (n - 2)) t^7 + (-15 + 9 (n - 2)) t^8+ (-9 + 6 (n - 2)) t^9 \\
&& + (-4 + 3 (n - 2)) t^{10} + (-1 + (n - 2)) t^{11}\\
&=&(1 + t)^2 (-1 - 3 t + n t) (1 + t^2) (1 - t + t^2) (1 + t + t^2)^2,
\end{eqnarray*}
\begin{eqnarray*}
H_2(t) & = & -1 + (-4 + (n - 2)) t + (-9 + 5 (n - 2)) t^2 + (-15 + 14 (n - 2)) t^3 \\
&&+ (-20 + 29 (n - 2)) t^4 + (-23 + 49 (n - 2)) t^5 + (-23 + 66 (n - 2)) t^6 \\
&&+ (-20 + 71 (n - 2)) t^7 + (-15 + 61 (n - 2)) t^8 + (-9 + 40 (n - 2)) t^9 \\
&&+ (-4 + 19 (n - 2)) t^{10} + (-1 + 5 (n - 2)) t^{11}\\
&=& (1 + t)^2 (1 + t^2) (1 + t + t^2) \\
&&(-1 - 3 t + n t - 5 t^2 + 2 n t^2 - 
   7 t^3 + 3 n t^3 - 9 t^4 + 4 n t^4 - 11 t^5 + 5 n t^5).
\end{eqnarray*}

Now we assume that $n \geqq 6$.
Then all coefficients of $H_1(t)$ and $H_2(t)$ except their constant terms
are non-negative so that they have unique zeros in $(0, \infty)$.
The following inequalities
$$
H_1(\frac{1}{n-3})=0, \;\; H_2(\frac{1}{n-1})=-\frac{6}{(n-1)^5}<0
$$
guarantee that the zero of $H(t)$ is located in $[\frac{1}{n-1}, \frac{1}{n-3}]$.
Combining with the similar result for $n=4,5$ in the previous section, we have the following theorem
which is our main result.

\begin{theo}
The growth rate $\tau$ of an ideal Coxeter polytope with $n$ facets in $\mathbb{H}^3$ satisfies
\begin{equation} \label{eq: main}
n-3 \leqq \tau \leqq n-1.
\end{equation}
\end{theo}


\begin{cor}
An ideal Coxeter polytope $P$ with $n$ facets in $\mathbb{H}^3$ is right-angled if and only if
its growth rate $\tau$ is equal to $n-3$.
\end{cor}
\proof
The factor $H(t)$ of the denominator polynomial $G(t)=(t-1)H(t)$ of the growth function of $P$ is equal to
$H_1(t)$ if and only if $b=c=d=0$, which means that all dihedral angles are $\pi/2$.
\qed

From the inequality (\ref{eq: main}), 
we see that the growth rate $\tau$  of an ideal Coxeter polytope with $n$ facets with $n \geqq 6$ 
satisfies $\tau \geqq 3$.
Therefore combining with the result of growth rates for $n=4,5$ shown in the previous section,
we also have the following corollary (see also \cite{N}, Theorem 4).

\begin{cor} 
The minimum of the growth rates of 3-dimensional hyperbolic ideal Coxeter polytopes is 
 $0.492432^{-1}=2.03074$, which is uniquely realized by the ideal simplex with 
$p=q=s=2$.
\end{cor}

\section{Acknowledgement}
The authors thank Jun Nonaka for explaining his paper \cite{N}.


\end{document}